\documentclass[titlepage,12pt]{article}
\usepackage{amssymb}
\usepackage{latexsym}
\newcommand{\ben}{\begin{enumerate}}
\newcommand{\een}{\end{enumerate}}
\newcommand{\ble}{\begin{lem}}
\newcommand{\ele}{\end{lem}}
\newcommand{\bth}{\begin{thm}}
\renewcommand{\eth}{\end{thm}}
\newcommand{\bpr}{\begin{prop}}
\newcommand{\epr}{\end{prop}}
\newcommand{\bco}{\begin{cor}}
\newcommand{\eco}{\end{cor}}
\newcommand{\bcon}{\begin{conj}}
\newcommand{\econ}{\end{conj}}
\newcommand{\bde}{\begin{defn}}
\newcommand{\ede}{\end{defn}}
\newcommand{\bex}{\begin{exa}}
\newcommand{\eex}{\end{exa}}
\newcommand{\barr}{\begin{array}}
\newcommand{\earr}{\end{array}}
\newcommand{\btab}{\begin{tabular}}
\newcommand{\etab}{\end{tabular}}
\newcommand{\beq}{\begin{equation}}
\newcommand{\eeq}{\end{equation}}
\newcommand{\bea}{\begin{eqnarray*}}
\newcommand{\eea}{\end{eqnarray*}}
\newcommand{\bce}{\begin{center}}
\newcommand{\ece}{\end{center}}
\newcommand{\bpi}{\begin{picture}}
\newcommand{\epi}{\end{picture}}
\newcommand{\bfi}{\begin{figure} \begin{center}}
\newcommand{\efi}{\end{center} \end{figure}}

\newcommand{\bsl}{\begin{slide}{}}
\newcommand{\esl}{\end{slide}}

\newcommand{\bib}{thebibliography}
\newcommand{\pf}{{\bf Proof}\hspace{7pt}}

\newcommand{\qqed}{\qquad\rule{1ex}{1ex}}
\newcommand{\Qqed}{\qquad\rule{1ex}{1ex}\medskip}

\newcommand{\hso}[1]{\hspace{-1pt}}

\newcommand{\qmq}[1]{\quad\mbox{#1}\quad}

\newcommand{\emp}{\emptyset}

\newcommand{\sbe}{\subseteq}

\newcommand{\spe}{\supseteq}

\newcommand{\iso}{\cong}

\newcommand{\zh}{\hat{0}}
\newcommand{\oh}{\hat{1}}

\newcommand{\lte}{\unlhd}

\newcommand{\case}[4]{\left\{\barr{ll}#1&\mbox{#2}\\#3&\mbox{#4}\earr\right.}

\def\<{\langle}
\def\>{\rangle}

\newcommand{\ree}[1]{(\ref{#1})}

\newcommand{\ra}{\rightarrow}

\newcommand{\mut}{\tilde{\mu}}

\newcommand{\De}{\Delta}

\newcommand{\bbZ}{{\mathbb Z}}
\newcommand{\cA}{{\cal A}}

\newcommand{\cB}{{\cal B}}

\newcommand{\cC}{{\cal C}}

\newcommand{\cF}{{\cal F}}

\newcommand{\cI}{{\cal I}}

\newcommand{\cJ}{{\cal J}}

\newcommand{\Ht}{\tilde{H}}

\newcommand{\rk}{\mathop{\rm rk}\nolimits}

\newcommand{\cho}{\choose}

\newcommand{\ds}{\displaystyle}

\newcommand{\act}{\mathop{\rm Act}\nolimits}
\newcommand{\Act}{\mathop{\rm Act}\nolimits}

\newcommand{\atoms}{{\cal A}}

\newcommand{\bases}{{\cal B}}

\newcommand{\bigjoin}{\bigvee}

\newcommand{\dfn}{\em}

\newcommand{\ext}{\mathop{\rm Ext}\nolimits}
\newcommand{\Ext}{\mathop{\rm Ext}\nolimits}

\newcommand{\join}{\vee}

\newcommand{\lee}{\le^{\rm ext}}
\newcommand{\lei}{\le^{\rm int}}

\newcommand{\mb}{\mathop{\rm MaxBas}\nolimits}

\newcommand{\maxelt}{\hat{1}}
\newcommand{\minelt}{\hat{0}}
\newcommand{\nerve}{{\rm Nerv}}

\newcommand{\redhom}{{\mathop{\rm \tilde H}}}

\newcommand{\redeulch}{{\mathop{\tilde \chi}}}

\newcommand{\AI}[1]{\item[\rm{(#1)}]} 
\newcommand{\NI}[1]{\item[$#1$.]}
\newcommand{\IN}{{I\! N}}

\newtheorem{thm}{Theorem}[section]
\newtheorem{prop}[thm]{Proposition}
\newtheorem{cor}[thm]{Corollary}
\newtheorem{lem}[thm]{Lemma}
\newtheorem{conj}[thm]{Conjecture}
\newtheorem{exa}[thm]{Example}

\begin{document}
\pagestyle{empty}
\title{Topological properties of active orders for matroid bases}
\author{Rieuwert J. Blok \\
Department of Mathematics\\ 
Michigan State University\\
East Lansing, MI 48824-1027
USA\\
{\tt blokr@member.ams.org}\\[5pt]
and\\[5pt]
Bruce E. Sagan\\
Department of Mathematics\\ 
Michigan State University\\
East Lansing, MI 48824-1027\\
USA\\
{\tt sagan@math.msu.edu}}

\date{\today \\[1in]
	\begin{flushleft}
	\end{flushleft}
       }
\maketitle

\begin{flushleft} Proposed running head: \end{flushleft}
	\begin{center} 
Topology of active orders
	\end{center}

Send proofs to:
\begin{center}
Bruce E. Sagan \\ Department of Mathematics \\ Michigan State
University \\ East Lansing, MI 48824-1027\\[5pt]
Tel.: 517-355-8329\\
FAX: 517-432-1562\\
Email: sagan@math.msu.edu
\end{center}

\begin{abstract}
Las Vergnas~\cite{lv:aom} introduced several lattice structures on the
bases of an ordered matroid $M$ by using  their external and internal
activities.  He also noted~\cite{lv:pc} that when computing the
M\"obius function of these lattices, it was often zero, although he
had no explanation for that fact.  The purpose of this paper is to
provide a topological reason for this phenomenon.  In particular, we
show that the order complex of the external lattice $L$ of $M$ is
homotopic to the independence complex of the restriction $M^*|T$ where
$M^*$ is the dual of $M$ and $T$ is the top element of $L$.  We then
compute some examples showing that this latter complex is often 
contractible which forces all its homology groups, and thus its M\"obius
function, to vanish.
A theorem of Bj\"orner~\cite{Bj1992} also helps us to calculate the
homology of the matroid complex.
\end{abstract}

\pagestyle{plain}
\section{The external and internal orders}

In September of 2001, there was a conference on Tutte Polynomials and
Related Topics at the Centre de Recerca Matematica in Barcelona,
Spain.  At the meeting, Michel Las Vergnas gave a talk about three
lattice structures which he had imposed on the bases of an ordered matroid
using external and internal activity~\cite{lv:aom}.  During the
question and answer period that followed, one of us (Sagan), asked if
Las Vergnas knew anything about the M\"obius function of these
lattices.  Las Vergnas replied that he had computed some examples
and noted that the value was often zero, but did not have an
explanation for that fact.

In this paper, we will give a topological reason for Las Vergnas'
observation.  The rest of this section will be devoted to developing the
definition and some basic properties of the external lattice, $L$.
In the next section, we derive some results about the structure of $L$
which will be useful in working with its order complex $\De$.  In
particular, we give a simpler formula for the join operator than was
given by Las Vergnas.  The third section contains our main theorem,
showing that $\De$ is
homotopic to the independence  complex $\IN$ of the restriction
$M^*|T$ where $M^*$ is the dual of $M$ and $T$ is the top element of
$L$.  In section 4, we compute some examples showing that $\IN$ is
often contractible which forces all its homology groups, and thus its M\"obius
function, to be zero.  A characterization of the homology of $\IN$ due
to Bj\"orner~\cite{Bj1992} is recalled in the next section and used
for the calculation of yet more examples.  The final section
contains a couple of open problems.

Let $M$ be a matroid on a finite set $E$.
We denote the bases and independent sets of $M$ by  $\cB=\cB(M)$
and $\cI=\cI(M)$, respectively. 
We say that $M$ is {\dfn ordered} if $E$ is linearly ordered.
From now on all matroids will be ordered.

Given a set $F\sbe E$ we say that $e\in E$ is {\dfn active with respect to
 $F$} if there is a circuit $C(F;e)\sbe F\cup\{e\}$ in which $e$ is
minimal with respect to the ordering on $E$.
Let 
$$
\act_M(F)=\{e\ :\ \mbox{ $e$ is active with respect to $F$}\}.
$$
Note that we include the possibility that $e\in F$.  Note also
that we will often write one-element sets without the set braces and
drop $M$ as a subscript if the matroid is clear from context. 

For $F\sbe E$ we define
$$
\ext_M(F)=\act_M(F)-F.
$$
The elements of $\ext_M(F)$ are called {\dfn externally active with respect
 to $F$}.
This coincides with the usual notion of externally active elements with respect
to an element of $\bases$.  

Las Vergnas defined the external lattice of $M$ in a
manner equivalent to the following.  For $A,B\in\bases$, define
$$
\mbox{$A\lee_M B$ if and only if $A\sbe B\cup \Ext_M(B)$.}
$$
It was proven in~\cite{lv:aom} that, when augmented with a minimum
element $\zh$, the resulting order is in fact 
a graded lattice with rank function 
\beq
\label{rho}
\rho_M(B)=|\Ext_M(B)|+1.
\eeq
We will denote this lattice by $L(M)$ or simply $L$.  It is important
to remember that, even though our notation does not show it, this
lattice structure depends not the ordering of the base set of $M$.

Let us give two simple examples of these
lattices by using the cycle matroid of a multigraph $G=(V,E)$.  
Because of symmetry, these particular lattices do not depend on the
ordering of the edge set.  If $G$ is
the circuit on $n$ vertices, then $L$ consists of an antichain of
$n-1$ elements with a minimum and maximum adjoined.  At the other
extreme, if $G$ consists of two vertices with $n$ edges between them then
$L$ is a chain of $n+1$ elements.

Returning to our general exposition, let $M^*$ be the dual matroid of $M$.
We turn $M^*$ into an ordered matroid using the order already given on $E$.
Las Vergnas~\cite{lv:aom} also defined another ordering 
$\lei_M$ on $\bases(M)$ by
\beq
\label{int}
A\lei_M B\iff (E-B)\lee_{M^*} (E-A).
\eeq
We should note that one can also define $\lei_M$ using the internal
activity of bases of $M$ (which also eliminates the need to pass to
$M^*$), but~\ree{int} will be more convenient for our purpose.
When augmented with a maximum element $\oh$, the resulting order 
is called the {\em internal order}.  Directly from the
definitions, we see that this structure is just the order-theoretic
dual of $L(M^*)$.  Since the dual of a lattice has the
same homology as the original lattice, we will restrict ourselves to
external orders.  For that reason, we will also drop the ${\rm ext}$
superscript.

It will be useful in the sequel to have the following characterization,
due to Las Vergnas~\cite[Proposition 3.1]{lv:aom}. of the external order.
\begin{prop}[Las Vergnas]  
\label{lasvergnas}
Let $A$, $B$ be two bases of an ordered matroid $M$. 
Then $A\leq B$ if and only if  $B$ is the lexicographically maximum 
base of $M$ contained in $A\cup B$ (where elements of a base are
listed in increasing order).\qqed
\end{prop}

In the aforementioned paper it was shown that the 
 number of elements at a given rank in $L(M)$ does not depend on
the particular order on $E$, but that the lattice itself does.
We wish to give some measure of how $L(M)$ depends on the order on $E$.

\bpr\label{notexact}
Let $\lte$ and $\lte'$ be linear orders on $E$.  Given a matroid on
$E$, let $M$ and $M'$ be the corresponding ordered matroids.
Suppose that $\act(M)=\act(M')$ and that $\lte$, $\lte'$ when restricted to
this set are same.
Then 
$$
L(M)\cong L(M').
$$
\epr
\pf\
We prove that the identity map from $\cB(M)$ to $\cB(M')$ induces
a lattice isomorphism of $L(M)$ with $L(M')$. 
So we need to show that  for  $A,B\in\bases(M)=\bases(M')$  we have
$A\sbe B\cup \Ext_M(B)$ if and only if $A\sbe B\cup \Ext_{M'}(B)$.
Clearly it suffices to  have $\Ext_M(B)=\Ext_{M'}(B)$.  We will show 
$\Ext_M(B)\sbe\Ext_{M'}(B)$ and then the reverse inclusion follows by
symmetry.  Now take $a\in\Ext_M(B)$ and let $C$ be the unique cycle in
$B\cup a$.  So $a$ is the $\lte$-minimum in $C$ and it suffices to
show that it is also the $\lte'$-minimum.  Let $a'$ be this $\lte'$-minimum.
Then 
$a,a'\in\Act(M)=\Act(M')$ with $a\lte a'$ and $a'\lte' a$.  Since the
two orderings agree on this set, $a=a'$ and we are done.\qqed

\section{Sublattices and the join operator}

Fix a subset $F\sbe E$ and let $K=M|F$ be the restriction of $M$ to $F$.
Note that it is an ordered matroid with respect to the ordering
induced on $F$ by $E$.  We will say that $K$ is {\em spanning\/} if $F$ is a
spanning set of $M$, that is, $F$ contains a base of $M$.
We will show that the lattice for a spanning matroid is closely
related to that of the parent matroid.  But first we need a lemma.

\ble
\label{ext}
Suppose that $F\sbe E$ and $K=M|F$.  Then for any $J\sbe F$ we have
\begin{itemize}
\AI{a}
$\ds \act_K(J)=\act_M(J)\cap F$, and as a consequence
\AI{b}
$\ds \Ext_K(J)=\Ext_M(J)\cap F$. 
\end{itemize}
\ele
\pf
(a)  The fact that $\act_K(J)\sbe \act_M(J)\cap F$ is clear from the
definitions.  For the opposite inclusion, suppose
$e\in\act_M(J)\cap F$.  Then there is a circuit 
$C\sbe J\cup e$ in which $e$ is minimal. 
But then $C\sbe F$ and $e$ is minimal with respect to the 
ordering induced on $F$ so that $e\in\act_K(J)$.

Part (b) follows immediately from part (a).
\qqed

\begin{cor}\label{cor:sublattice} 
Suppose that $K=M|F$ is spanning.
Then the inclusion $\bases(K)\sbe\bases(M)$ induces an inclusion
 $$L(K)\sbe L(M).$$
\end{cor}
\pf
Suppose $A,B\in\bases(K)$.   We prove that $A\leq_M B$ if and only if
 $A\leq_K B$.
By definition, $A\leq_M B$ if and only if $A\sbe B\cup\ext_M(B)$.
Since $A,B\sbe F$ this happens if and only if 
$A\sbe B\cup(\ext_M(B)\cap F)$.
By the previous lemma,
 $B\cup(\ext_M(B)\cap F)=B\cup\ext_K(B)$.
So we are done.\qqed

Following Las Vergnas~\cite{lv:aom}, for a spanning subset $A\sbe E$ we define 
$$
\mb{A}=A-\act(A). 
$$
Alternatively, one can define this as the lexicographically maximum base
of $M$ contained in $A$, using the convention of
Proposition~\ref{lasvergnas}. 
We obtain the maximum element of $L=L(M)$ as 
$$
T=\mb{E} 
$$
and reserve the notation $T$ for this top element.
Las Vergnas gave a formula for the join operator $\join$ for two elements
of $L$ using the $\mb{}$ operator.
Using Corollary~\ref{cor:sublattice} we give a slight but useful
simplification of his result, at the same time extending it to the
join of an arbitrary  number of elements in $L$.

\begin{cor}\label{cor:joinismaxbasofunion}
The join of elements $B_i\in\bases(M)$ ($i=1,2,\ldots,m$) in $L(M)$
is given by 
$$\bigjoin_{i=1}^m B_i=\mb\left(\bigcup_{i=1}^m B_i\right)$$
\end{cor}
\pf
Let $K=M|F$ where $F=\bigcup_{i=1}^m B_i$ and let $S=\mb(F)$.
We must prove that $S=\bigjoin_{i=1}^m B_i$.
First of all, for all $i$ we have $B_i\leq_K S$ because 
 $S$ is the maximal element of $L(K)$.
By Corollary~\ref{cor:sublattice} this means $B_i\leq_M S$
 for all $i$.

Now suppose $T\in\bases(M)$ satisfies
 $B_i\leq_M T$ for all $i$. Then $B_i\sbe T\cup\ext_M(T)$  so that
 $F=\bigcup_{i=1}^m B_i\sbe T\cup \ext_M(T)$.
But $S\sbe F\sbe T\cup\ext_M(T)$ and so by we have $S\leq_M T$.
Thus $S=\bigjoin_{i=1}^m B_i$.\Qqed

We denote the set of atoms of $L(M)$ by $\atoms=\cA(M)$. By~\ree{rho},
these are precisely the bases $B$ for $M$ with $\ext(B)=\emptyset$.

\begin{cor}
\label{cor:baseformaxbas}
Let $\cA'\sbe \atoms$.
Then $\bigjoin_{B\in \cA'}B=T$ if and only if every element of $T$ is contained
 in some element $B\in\cA'$.
\end{cor}
\pf
This follows from Corollary~\ref{cor:joinismaxbasofunion} and the following
 observation which is needed for the ``if'' direction.
Suppose $T\sbe F$ for some $F\sbe E$. 
Then since $T\cap\act(F)\sbe T\cap\act(E)=\emptyset$ we have
   $T\sbe \mb(F)$.
Also, if $F$ is spanning, then $\mb(F)$ is a base for $M$.
Since $T$ is also a base for the matroid $M$, we find $T=\mb(F)$.
\qqed

\medskip

The inclusion in Corollary~\ref{cor:sublattice}  does not preserve
the rank function in general.  But it does under certain circumstances.

\begin{lem}
\label{lem:inclusion preserves rank}
If $K=M|F$ is spanning and $B\in\bases(K)\sbe\bases(M)$,
then the following hold.
\begin{itemize}
\AI{a} We have $\rho_{K}(B)=\rho_{M}(B)$ if and only if $\ext_M(B)\sbe F$.
\AI{b} If $F\spe E-T$, then the inclusion $L(K)\sbe L(M)$ 
 preserves rank.
\AI{c} If $f<e$ for all $f\in F$ and $e\in E-F$,
 then the inclusion $L(K)\sbe L(M)$ preserves rank.
\end{itemize}
\end{lem}
\pf
(a) We have $\rho_{K}(B)=|\ext_{K}(B)|+1$ and
 $\rho_{M}(B)=|\ext_{M}(B)|+1$. Now Lemma~\ref{ext} completes the proof.

(b) This follows from part (a) since for any $A\sbe E$ we have
$\ext_M(A)\sbe \act_M(E)=E-T$.

(c) This also follows from part (a) since the assumption implies that no element
 of $E-F$ can be externally active with respect to any subset of $F$.
\qqed

Given a subset $F\sbe E$ and an ordering on $F$ we can always
 define an ordering on $E$ such that the
 condition in (c) of Lemma~\ref{lem:inclusion preserves rank} holds. 
Thus we have proved the following observation.

\begin{cor}\label{cor:supermatroid}
Let $K$ be an ordered matroid on a set $F$.
If $M$ is an unordered matroid on a set $E\spe F$ such that 
 $K=M|F$ and $K$ is spanning, then we can find an ordering on $E$
inducing a rank-preserving
 inclusion $L(K)\sbe L(M)$.
\end{cor}

In particular if $K$ is the cycle matroid of a connected graph $H$
with edge set $F$,  then for $M$ we can take the cycle matroid of the
complete graph on  the vertex set of $H$.

\section{The homotopy equivalence}
\label{h}

In this section we study the reduced homology of the order complex of the 
lattice $L(M)$.
We will show that there is a homotopy equivalence between the order complex
of $L(M)$ and the  independence complex of $M^*$
restricted to $T$.  This will we used in the next section to explain
Las Vergnas' observation about the M\"obius function of $L(M)$. 

Let $L$ be a finite lattice with minimum and maximum elements $\zh$ and $\oh$,
respectively.  We denote by $\De(L)$, or simply $\De$, the 
{\dfn order complex} of $L$, that is, 
the abstract simplicial complex on the set $L-\{\zh,\oh\}$ whose 
faces are the nonempty chains in $L-\{\minelt,\maxelt\}$ 
ordered by inclusion.  If $L=L(M)$ for some matroid, then we will also
use the notation $\De(M)=\De(L(M))$.

There is another abstract simplicial complex associated with a matroid.
The {\dfn independence complex} of $M$, denoted $\IN(M)$, is the
simplicial complex  of nonempty independent subsets of $M$.
Let $T'$ be the elements of $T$ that are independent as singleton sets
in $M^*$.  Then $\IN(M^*|T)=\IN(M^*|T')$.
Note that the elements $e\in E$ which are not independent in $M^*$ are
precisely those which are contained in every base for $M$.  Our main 
theorem relates the two complexes we have defined.
In it, $\redhom_i(\De)$ will denote the reduced $i$-dimensional
homology group of  a complex $\De$ with coefficients in $\bbZ$  (see e.g.\
Stanley~\cite[Ch.3]{Sta1986}).

\begin{thm}
\label{De=IN}
We have a homotopy equivalence
$$
\De(M)\simeq\IN(M^*|T).
$$ 
So, for all $i\ge-1$, we have an isomorphism in homology
$$\redhom_i(\De(M))\cong\redhom_i(\IN(M^*|T)).$$
\end{thm}
Note that this result implies that the homotopy type of the
order complex depends only on the maximum base $T$.
We will prove  Theorem~\ref{De=IN} using the next two propositions.

Let $L$ be an arbitrary lattice with atom set $\cA$.  Let $\cJ=\cJ(L)$ be the
abstract simplicial complex of all subsets of $\cA$ whose join is not
$\oh$.  The following is a  theorem of
Lakser~\cite{lak:hl} later generalized by Bj\"orner~\cite{bjo:htp} and
Segev~\cite{seg:ocp}. 
\begin{prop}
\label{De=N}
For any lattice $L$
$$
\De(L)\simeq\cJ(L).\qqed
$$
\end{prop}

Let $\cF$ be an abstract simplicial complex on a finite set $F$.
A {\em facet covering} of $\cF$ is a multiset of facets
 $\cC=\{F_0,F_1,\ldots,F_n\}$
 such that every face of $\cF$ is contained in some $F_i$.
The {\em nerve} $\nerve(\cC)$ of the covering is the simplicial complex
 on the vertex set $I=\{0,1,2,\ldots,n\}$ where
 a subset $J\sbe I$ is a face if and only if $\bigcap_{j\in J} F_j$ is a
 face of $\cF$.
As  will be seen, the nerve of a certain covering of $\cJ(L)$ is isomorphic
 to $\IN(M^*|T)$.

But first we must show that  $\cF$ and $\nerve(\cC)$ are 
 the same up to homotopy.  
Note that every  nonempty intersection of facets of $\cF$
 is again a face of $\cF$.
Thus the intersections $\bigcap_{j\in J}F_j$
 are contractible as subspaces of $\cF$ and hence are acyclic.
Thus the hypotheses of the Nerve Theorem of Borsuk and Folkman are
satisfied (see (10.6) in Bj\"orner~\cite{bjo:tm}) and we obtain our second
proposition. 
\begin{prop}
\label{prop:intersection complex}
Let $\cF$ be a simplicial complex on a set $F$ and let $\cC$ be a 
 facet covering.
Then
$$
\cF\simeq\nerve(\cC).\qqed
$$
\end{prop}
\pf {\bf (of Theorem~\ref{De=IN})}
Combining Propositions~\ref{De=N} and~\ref{prop:intersection complex}
for any facet covering $\cC$ of $\cJ$ we have
$$
\De\simeq\cJ\simeq\nerve(\cC).
$$  
So it suffices to show that we can find a facet covering $\cC$ such that
 $\nerve(\cC)$ and $\IN(M^*|T)$ are isomorphic as simplicial complexes.

Suppose $T'=\{t_0,t_1,\ldots,t_n\}$ and recall that
$\IN(M^*|T)=\IN(M^*|T')$. 
For $0\le i\le n$, define $F_i=\{A\in\atoms\ :\ A\sbe E-\{t_i\}\}$.
Then it follows from Corollary~\ref{cor:baseformaxbas} that these are
 the facets of $\cJ$, possibly with repetitions.  
Let $\cC$ be the corresponding facet covering of $\cJ$.
We can now define a bijection
$\phi:\IN(M^*|T')\ra\nerve(\cC)$ as follows.  If $S\sbe T'$ then let
$$
\phi(S)=J=\{j\ :\ t_j\in S\}.
$$
Clearly $\phi$ is  a bijection between subsets of  $T'$ and subsets of $I$.
We claim that $\phi$ restricts to a well-defined isomorphism between
the respective complexes, that is,
$\bigcap_{j\in\phi(S)}F_j\neq\emp$  if and only if $S$ is independent
in $M^*|T'$.  This is because $S$ is independent in $M^*|T'$ if and
only if  $E-S$ contains  a base for $M$ which, by
Lemma~\ref{lem:inclusion preserves rank}(b), is  
 equivalent to $E-S$ containing an atom for $L$.
This completes the proof of the isomorphism and of Theorem~\ref{De=IN}.
\qqed

\section{Applications}

We are now ready to explain the empirical observation of Las Vergnas
that the M\"{o}bius function $\mu$ of the external lattice $L(M)$
often satisfies $\mu(L(M))=0$.  It is known that, given any finite lattice
$L$ with minimum element $\zh$, maximum element $\oh$, and M\"{o}bius
function $\mu$, one has  
\beq
\label{chi}
\mu(L):=\mu_L(\zh,\oh)=\redeulch(\De)
=\sum_{i=-1}^\infty(-1)^i\dim \redhom_i(\De) 
\eeq
where $\De$ is the order complex of $L$ and 
$\redeulch$ is the reduced Euler characteristic.
This equation together with Theorem~\ref{De=IN}
can be used to show that a number of external activity lattices have
M\"obius function zero.  We will use the notation $\Ht_i(M)$ and
$\mu(M)$ for $\Ht_i(\De(M))$ and $\mu(L(M))$, respectively.  We will
also use $\rk(M)$ for the rank of the matroid $M$.  This should not be
confused with the rank function $\rho$ for the lattice $L(M)$.

\bpr
\label{ball-sphere}
Let $M$ be an ordered matroid with maximum base $T$ and rank
$r=\rk(M)\ge1$.
\ben
\AI{a} Suppose that $M|(E-T)$ is spanning. Then
$$
\mbox{$\Ht_i(M)=\{0\}$ for all $i\ge-1$}\qmq{and} \mu(M)=0.
$$
\AI{b} Suppose that $M|(E-S)$ is spanning for all proper
subsets $S\subset T$ but is not spanning for $S=T$.  Then
$$
\Ht_i(M)=\case{\bbZ}{if $i=r-2$,}{\{0\}}{else,}\qmq{and} 
\mu(M)=(-1)^{r-2}.
$$
\een
\epr
\pf
Under the first (respectively, second) hypothesis, $\IN(M^*|T)$ is
homologically an $(r-1)$-ball (respectively, $(r-2)$-sphere).  The
conclusions now follow from Theorem~\ref{De=IN} and equation~\ree{chi}.\qqed

\medskip

As an example, consider the cycle matroid of a graph $G$ where, as
usual, the edge set $E=E(G)$ has been linearly ordered.  In this
case we will use $G$ in our notation everywhere we used $M$ before.
If $H$ is a subgraph of $G$, then we will call a vertex $v$ of $H$
{\it internal\/} if every edge of $G$ containing $v$ is present in
$H$.
\bco
Let $K_n$ be an ordered complete graph on $n$ vertices, $n\ge2$, and
let $T$ be its lexicographically maximal spanning tree.  
\ben
\AI{a} If $T$ has no internal vertex then
$$
\mbox{$\Ht_i(K_n)=\{0\}$ for all $i\ge-1$}\qmq{and} \mu(K_n)=0.
$$
\AI{b} If $T$ has an internal vertex then
$$
\Ht_i(K_n)=\case{\bbZ}{if $i=n-3$,}{\{0\}}{else,}\qmq{and} 
\mu(K_n)=(-1)^{n-3}.
$$
\een
\eco
\pf
If $T$ has no internal vertex, then $K_n-E(T)$ is connected and the
 hypotheses of Proposition~\ref{ball-sphere}~(a) are satisfied.  
If $T$ has an internal vertex, then $K_n-E(S)$ is connected for all $S\sbe T$, 
 except for $S=T$.
Thus the hypotheses of Proposition~\ref{ball-sphere}~(b)
 are fulfilled.\qqed 

\medskip

As a result of this corollary, we can see that $\De(M)$ is not, in general,
shellable (even though $\IN(M^*|T)$ always is, see
Bj\"orner~\cite[Theorem 7.3.3]{Bj1992}).
If $\De$ is a shellable simplicial complex pure
of dimension $d$, then $\De$ is topologically a wedge of $d$-spheres and
so only has homology in dimension $d$.  So if a finite lattice $L$ graded of
rank $\rho$ is shellable, then it only has homology in dimension $\rho-2$
(since we remove $\zh$ and $\oh$).  But in $L(M)$ we have
$$
\rho(L(M))=\rho(T)=|\Ext(T)|+1=|E-T|+1.
$$
In particular
$$
\rho(L(K_n))={n\choose2}-(n-1)+1={n-1\choose2}+1.
$$
But from the previous corollary, if $T$ has an internal vertex then
$L(K_n)$ has homology in dimension $n-3<{n-1\choose2}-1$ for $n\ge4$.

\medskip

Here is another family of matroids that have zero M\"obius function.
\bpr
\label{cone}
Let $M$ be an ordered matroid with maximum base $T$ and suppose there
 is $t\in T$ such that $\rk(E-T)=\rk((E-T)\cup t)$.  Then
$$
\mbox{$\Ht_i(M)=\{0\}$ for all $i\ge-1$}\qmq{and} \mu(M)=0.
$$
\epr
\pf
Suppose that $t\in T$ satisfies $\rk(E-T)=\rk((E-T)\cup t)$.  
This means that if a base $B\in\bases(M)$ intersects $T$ minimally,
 then $t\not\in B$.
That is, $t$ is not contained in any base of the contraction $M.T$ and
hence is contained  in every base of $M^*|T$.
Thus $\IN(M^*|T)$ is a cone with vertex $t$.
The result follows.\qqed

\medskip

For application in our examples, note that for the cycle matroid of a
graph $G$, the hypothesis of Corollary~\ref{cone} just says that the
edge $t\in T$ connects two vertices in the same component of $G-E(T)$.
We first consider the {\it $n$-fan}, $F_n$, which is obtained from a
path with $n$ vertices by adding an additional vertex adjacent to
every vertex of the path.  More explicitly, $F_n=(V,E)$ where 
$V=\{0,1,\ldots,n\}$ and 
$$
E=\{01,02,\ldots,0n\}\uplus\{12,23,\ldots,(n-1)n\}.
$$
We always write our edges with the smaller vertex
first and order them lexicographically.  Then  
$$
E(T)=\{0n,12,23,\ldots,(n-1)n\}.
$$
It is easy to see that if $n\ge3$ then the edge $t=12$ satisfies the
component criterion of the first sentence in this paragraph.

Next consider the {\it $n$-triangle graph}, $T_n$, gotten by gluing
together $n$ copies of $K_3$ along a common edge.  To set notation,
let
$$
E=\{e_0,e_1,\ldots,e_{2n}\}
$$
where the $i$th triangle has edges $\{e_0,e_i,e_{n+i}\}$ and edges are
ordered by their subscripts.  Now
$$
T=\{e_n,e_{n+1},\ldots,e_{2n}\}
$$
So if $n\ge3$ then the edge $t=e_{n+1}$ will satisfy
the component criterion.  By Proposition~\ref{cone}, we have proved
the following.

\bpr
For the given orderings and $n\ge3$ we have
$$
\mbox{$\Ht_i(F_n)=\Ht_i(T_n)=\{0\}$ for all $i\ge-1$}\qmq{and} 
\mu(F_n)=\mu(T_n)=0.\qqed
$$
\epr

\section{A theorem of Bj\"orner}

A theorem of Bj\"orner~\cite[Theorem 7.8.1]{Bj1992} characterizes the reduced
homology of $\IN(M)$ for any matroid $M$ and can be used in
conjunction with Theorem~\ref{De=IN} for computations.
To state it, we will  need the lattice of flats of $M$ 
which will be denoted $L_F=L_F(M)$ to distinguish it from the external
activity lattice.  Also, define the {\it reduced M\"obius
function\/} of $M$ to be
$$
\mut(M)=\case{|\mu(L_F(M))|}{if $M$ is loopless,}{0}{else.}
$$
\bth[Bj\"orner]
If $r=\rk(M)$ then
$$
\Ht_i(\IN(M))\iso\case{\bbZ^{\mut(M^*)}}{if $i=r-1$,}{\{0\}}{else.}\qqed
$$
\eth

Now if $F\sbe E$, consider  $M.F$, the
contraction of $M$ to $F$.  Our interest stems from
the fact that $(M^*|F)^*=M.F$.  An immediate corollary of the previous
theorem and Theorem~\ref{De=IN} is as follows.
\bth
\label{Ht by mut(M.T)}
If  $r^*=\rk(M^*|T)$ then
$$
\Ht_i(\De(M))\iso\case{\bbZ^{\mut(M.T)}}{if $i=r^*-1$,}{\{0\}}{else.\qqed}
$$
\eth

\bco
\label{mutcor}
If $r=\rk(M)$ and $r^*=\rk(M^*|T)$ then
$$
\mu(M)=\case{(-1)^{r-1}\mu(L_F(M.T))}{if $M.T$ is loopless,}{0}{else.}
$$
\eco
\pf
Viewing $\mu(M)$ as the reduced Euler characteristic of 
 $\De(M)$ and using Theorem~\ref{Ht by mut(M.T)} we find
 $\mu(M)=(-1)^{r^*-1}\mut(M.T)$. 
So if $M.T$ has loops then $\mu(M)=0$ by definition of  $\mut$.
Otherwise, since $M.T=(M^*|T)^*$ and $|T|=r$, the rank of 
 $M.T$ and hence of $L_F(M.T)$ is $r-r^*$.
As $L_F(M.T)$ is a geometric lattice, the sign of $\mu(L_F(M.T))$
 is $(-1)^{r-r^*}$ and  canceling appropriate powers of $-1$ gives the
 desired conclusion.
\qqed

\medskip

Let us apply these results to some examples.

\paragraph{The uniform matroid}
Consider the {\dfn uniform matroid} $U_{n,k}$ on the 
 $n$-set $E$ whose collection of bases is 
 $$\cB(U_{n,k})=\{I\sbe E\ :\ |I|=k\}.$$

The lattice of flats $L_F(U_{n,k})$ consists of the subsets of $E$
 of cardinality strictly less than $k$ together with $E$ itself, 
 ordered by inclusion.
Thus $L_F(U_{n,k})$ is obtained from the Boolean lattice $B_n$ on $E$
 by deleting all elements of rank $l\ge k$, except the top element.
We will call this poset the {\em truncated Boolean algebra} (see
 Zhang~\cite{Zha1994})
Using the fact that, for any two subsets $A\sbe B\sbe E$, 
 the M\"{o}bius function of $B_n$ satisfies
 $$\mu(A,B)=(-1)^{|B-A|},$$
we find that
$$\mu(L_F(U_{n,k}))
 =-\sum_{i=0}^{k-1}(-1)^i {n\cho i} =(-1)^{k}{n-1\cho k-1}.$$

Now let $M=U_{n,k}$ for some $n>0$, and order $E$ linearly.
The top element $T$ of $L$ is some $k$-subset of $E$.
One verifies that $M^*|T$ is the uniform matroid $U_{k,r^*}$, where
$r^*=\min\{k,n-k\}$, and that $M.T$ is the uniform matroid $U_{k,k-r^*}$.
 
Suppose $k\le n/2$.  Then $r^*=k$ 
and only the empty set is independent in $M.T$.
Hence $M.T$ has loops, $\mut(M.T)=0$, and we have
$\redhom_i(\De)=\{0\}$ for all $i$, and $\mu(L)=0$.

Suppose instead that $k>n/2$ so that  $r^*=n-k$.  Then $M.T$ has no
loops and combining our computation  of $\mu(L_F(U_{n,k}))$ with
Theorem~\ref{Ht by mut(M.T)} and Corollary~\ref{mutcor} we have the
following result.  In it, we assume that ${j\choose i}=0$ if $i<0$.
\bpr
\label{Unk}
For any ordering of the uniform matroid $U_{n,k}$ we have
$$
\mbox{$\dim\Ht_i(U_{n,k})={k-1\choose 2k-n-1}$ if $i=n-k-1$ \quad and
\quad $\mu(U_{n,k})=(-1)^{n-k-1}{k-1\choose 2k-n-1}$.\qqed}
$$
\epr
Note that since $L(U_{n,k})$ has rank $n-k+1$, the complex
$\Delta(U_{n,k})$ is pure of  dimension $n-k-1$.  Apparently $\Delta(U_{n,k})$
only has homology in the top dimension. 

\paragraph{The wheel graph $W_n$}
Consider the {\dfn $n$-wheel graph}, $W_n$, obtained from an $n$-circuit
 $C$ by adding a vertex $v_0$ adjacent to all vertices of the circuit.
Let the edge set be ordered linearly and let $T$ be the top element of
$L(W_n)$. 

Suppose first that some edge $t\in T$ satisfies Proposition~\ref{cone} ,
i.e., $t$ connects two vertices in the same component of $W_n-E(T)$.
Then $\redhom_i(W_n)=\{0\}$ for all $i\ge-1$, and $\mu(W_n)=0$.

If there is no such edge, then $W_n-T$ is partitioned into connected
 components $C_0,C_1,\ldots,C_k$ as follows:
\begin{itemize}
\NI{1} $k=1$, $C_0=\{v_0\}$ and $C_1=C$, or 
\NI{2} $C_0$ is the union of triangles intersecting only in $v_0$,
           the components $C_1$, $C_2$,\ldots, $C_l$ are paths, possibly 
	   of length $0$, and every edge of $T$ meets $C_0$ and $C_i$ 
	   for some $i\ge 1$.
\end{itemize}

Let $T_i$ be the set of edges from $T$ joining $C_0$ to $C_i$.
Then by the above we have $T=\uplus_{i=1}^k T_i$.

Now $M.T$ is the cycle matroid of the graph with vertex set
 $\{C_0, C_1,\ldots, C_k\}$, where $T_i$ 
 represents a set of parallel edges joining the central vertex $C_0$
 to $C_i$.
Thus $M.T$ is the matroid of partial transversals
 of $T$ with respect to the family $\{T_i\}_{i=1}^k$.

We now determine $L_F(M.T)$.
The closed sets of $M.T$ are the unions of the sets $T_i$.
Thus $L_F(M.T)$ is the Boolean algebra $B_k$ on the set $\{T_i\}_{i=1}^k$.
Hence $\mu(L(M.T))=(-1)^k$.
Clearly $M^*|T=(M.T)^*$ has rank $n-k$ and so, 
using Theorem~\ref{Ht by mut(M.T)} and its corollary, we obtain the
following result
\bpr
Let $T$ be the top element of $L(W_n)$ for some ordering of the edges
of $W_n$.
\ben
\item  If there is an edge $t\in T$ satisfying Proposition~\ref{cone}
then
$$
\mbox{$\redhom_i(W_n)=\{0\}$ for all $i\ge-1$ \quad and \quad $\mu(W_n)=0$.}
$$
\item  If there is no such edge, then
$$
\dim \Ht_i(W_n)=\case{1}{if $i=n-k-1$,}{0}{else,}\qmq{and}
\mu(W_n)=(-1)^{n-k-1}.\qqed
$$
\een
\epr
Note that since $L$ has rank $n+1$, the complex $\Delta$ is pure of 
dimension $n-1$.  We have just shown that $\Delta$ has homology in
dimension $n-k-1$ and  since $k$ cannot be zero, this complex is not
shellable. 

\section{Open problems}

We observed  that the order complex for the uniform matroid has
homology in the correct dimension for it to be shellable.  It would be
nice to find an explicit shelling if one exists.  This would give a
way of re deriving Theorem~\ref{Unk}.

Forman~\cite{for:mtc} has introduced a discrete analogue of Morse
theory as a way of studying CW complexes by collapsing them onto
smaller, more tractable, complexes of critical cells.  These
techniques can be used to compute the homology of a complex even when
it is not shellable.  Are the non-shellable complexes which we have
considered amenable to Forman's technique?  

Las Vergnas~\cite{lv:aom} defined a third ordering on the bases of
an ordered matroid.  To state it, we first need one of his results.
\bpr[Las Vergnas]
If $A,B\in\cB(M)$ are distinct bases which are comparable in both
the external and internal orders, then either $A$ is smaller than $B$
in both  or $A$ is larger than $B$ in both.
\epr
Because of this proposition, we have a well-defined
{\it external-internal order $\le^{\rm exin}$\/}  on $\cB(M)$
given by 
$$
\mbox{$A\le^{\rm exin}_M B$ if and only if  $A\lee_M B$ or $A\lei_M B$}
$$
with corresponding lattice $L_{\rm inex}(M)$.
We have been unable to find an analogue of Theorem~\ref{De=IN} for
this lattice.  It would be very interesting to do so.

\noindent{\it Acknowledgment.}  We would like to thank Anders
Bj\"orner for helpful discussions and references.

\begin{\bib}{99}

\bibitem{bjo:htp} Anders Bj\"orner, Homotopy type of posets and
lattice complementation, {\it J. Combin.\ Theory Ser.\  A\/}
{\bf 30} (1981), 90--100. 

\bibitem{Bj1992} Anders Bj\"{o}rner, The Homology and Shellability of 
Matroids and Geometric Lattices in ``Matroid Applications,''
N. White ed., Encyclopedia of Mathematics and its Applications, Vol.\
40, Cambridge University Press, Cambridge, 1992, 226--283. 

\bibitem{bjo:tm} A. Bj\"orner, Topological methods, in ``Handbook of
Combinatorics,'' R. Graham, M. Gr\"otschel, and L. Lov\'asz eds.,
North-Holland, New York, NY, and MIT Press, Cambridge, MA, 1995,
1819--1872.

\bibitem{Bl2001} Andreas Blass,
Homotopy and homology of finite lattices, preprint (2001).

\bibitem{for:mtc} R. Forman, Morse theory for cell complexes, {\it
Adv.\ in Math.\/} {\bf 134} (1998), 90--145.

\bibitem{lak:hl} Harry Lakser, The homology of a lattice, {\it
Discrete Math.\/} {\bf 1} (1971--1972), 187--192.

\bibitem{lv:aom} M. Las Vergnas, Active orders for matroid bases,
{\it European J. Combin.\/} {\bf 22} (2001), 709--721.

\bibitem{lv:pc} M. Las Vergnas, personal communication.

\bibitem{BlSa1997} Andreas Blass and Bruce E. Sagan, M\"obius
functions of lattices, {\it Advances in Math.\/} {\bf 127} (1997),
94-123. 

\bibitem{seg:ocp}  Y. Segev, On the order complex of a prelattice,
{\it European J. Combin.\/} {\bf 18} (1997), 311--314.

\bibitem{Sta1986} Richard P. Stanley, Enumerative Combinatorics, Vol. I,
                  Cambridge Studies in Advanced Mathematics {\bf 49}, (1997)
                  Cambridge University Press, Cambridge.

\bibitem{We1976} D.J.A. Welsh, Matroid Theory, L.M.S. Monographs {\bf 8},
(1976)  Academic Press, London.

\bibitem{Wa1970} Andrew H. Wallace, Algebraic Topology: Homology and
Cohomology, (1970) W.A. Benjamin Inc., New York.

\bibitem{Zha1994} P. Zhang, Subposets of Boolean Algebras, PhD. thesis, 
                 Michigan State University (1994).
\end{\bib}
\end{document}